\documentclass{amsproc}%
\usepackage{graphicx}
\usepackage{amscd}
\usepackage{amsmath}
\usepackage{amsfonts}
\usepackage{amssymb}%
\theoremstyle{plain}

\numberwithin{equation}{section}

\begin{document}
\title[A power Cayley-Hamilton identity]{A power Cayley-Hamilton identity for $n\times n$ matrices over a Lie nilpotent
ring of index $k$}
\author{Jen\H{o} Szigeti}
\address{Institute of Mathematics, University of Miskolc \\
3515 Miskolc-Egyetemv\'{a}ros, Hungary}
\email{matjeno@uni-miskolc.hu}
\author{Szilvia Szil\'{a}gyi}
\address{Institute of Mathematics, University of Miskolc \\
3515 Miskolc-Egyetemv\'{a}ros, Hungary}
\email{matszisz@uni-miskolc.hu}
\author{Leon van Wyk*}
\address{Department of Mathematical Sciences, Stellenbosch University, Private Bag
X1,\hfill\break Matieland~7602, Stellenbosch, South Africa}
\email{LvW@sun.ac.za}
\thanks{The first named author was partially supported by the National Research,
Development and Innovation Office of Hungary (NKFIH) K119934.}
\thanks{*Corresponding author}

\begin{abstract}
For an $n\times n$ matrix $A$\ over a Lie nilpotent ring $R$ of index $k$,
with $k\geq2$, we prove that an invariant "power" Cayley-Hamilton identity
\[
\left(  I_{n}\lambda_{0}^{(2)}+A\lambda_{1}^{(2)}+\cdots+A^{n^{2}-1}%
\lambda_{n^{2}-1}^{(2)}+A^{n^{2}}\lambda_{n^{2}}^{(2)}\right)  ^{2^{k-2}}=0
\]
of degree $n^{2}2^{k-2}$ holds. The right coefficients $\lambda_{i}^{(2)}\in
R$, $0\leq i\leq n^{2}$ are not uniquely determined by $A$, and the cosets
$\lambda_{i}^{(2)}+D$, with $D$ the double commutator ideal $R[[R,R],R]R$
of~$R$, appear in the so-called second right characteristic polynomial
$p_{\overline{A},2}(x)$ of the natural image $\overline{A}$ of $A$ in the
$n\times n$ matrix ring~$\mathrm{M}_{n}(R/D)$ over the factor ring $R/D$:
\[
p_{\overline{A},2}(x)=(\lambda_{0}^{(2)}+D)+(\lambda_{1}^{(2)}+D)x+\cdots
+(\lambda_{n^{2}-1}^{(2)}+D)x^{n^{2}-1}+(\lambda_{n^{2}}^{(2)}+D)x^{n^{2}}.
\]

\end{abstract}
\maketitle

\noindent1. INTRODUCTION

\bigskip

The Cayley-Hamilton theorem and the corresponding trace identity play a
fundamental role in proving classical results about the polynomial and trace
identities of the $n\times n$ matrix algebra $\mathrm{M}_{n}(K)$ over a field
$K$ (see, for example, [Dr], [DrF] and [Row]).

In case of $\mathrm{char}(K)=0$, Kemer's pioneering work (see [K]) on the
T-ideals of associative algebras revealed the importance of the identities
satisfied by the $n\times n$ matrices over the Grassmann (exterior) algebra%
\[
E=K\left\langle v_{1},v_{2},...,v_{i},...\mid v_{i}v_{j}+v_{j}v_{i}=0\text{
for all }1\leq i\leq j\right\rangle
\]
generated by the infinite sequence of anticommutative indeterminates
$(v_{i})_{i\geq1}$.

Accordingly, the importance of matrices over non-commutative rings features
prominently in the theory of PI-rings; indeed, this fact has been obvious for
a long time in other branches of algebra, for example, in the structure theory
of semisimple rings. Thus any Cayley-Hamilton type identity for such matrices
seems to be of general interest.

In the general case (when $R$ is an arbitrary non-commutative ring with $1$)
Par\'{e} and Schelter proved (see [PaSch]) that a matrix $A\in\mathrm{M}%
_{n}(R)$ satisfies a monic identity in which the leading term is $A^{m}$ for
some large integer $m$, i.e., $m\geq2^{2^{n-1}}$. The other summands in the
identity are of the form $r_{0}Ar_{1}Ar_{2}\cdots r_{l-1}Ar_{l}$, with left
scalar coefficient $r_{0}\in R$, right scalar coefficient $r_{l}\in R$ and
"sandwich" scalar coefficients $r_{2},\ldots,r_{l-1}\in R$. An explicit monic
identity for $2\times2$ matrices arising from the argument of [PaSch] was
given by Robson in~[Rob]. Further results in this direction can be found
in~[Pe1] and~[Pe2].

Obviously, by imposing extra algebraic conditions on the base ring $R$, we can
expect \textquotedblleft stronger\textquotedblright\ identities in
$\mathrm{M}_{n}(R)$. A number of examples show that certain polynomial
identities satisfied by $R$ can lead to \textquotedblleft
canonical\textquotedblright\ constructions providing invariant Cayley-Hamilton
identities for $A$ of degree much lower than $2^{2^{n-1}}$.

If $R$ satisfies the polynomial identity%
\[
\lbrack\lbrack\lbrack\ldots\lbrack\lbrack x_{1},x_{2}],x_{3}],\ldots
],x_{k}],x_{k+1}]=0
\]
of Lie nilpotency of index $k$\ (with $[x,y]=xy-yx$), then for a matrix
$A\in\mathrm{M}_{n}(R)$, a left (and right) Cayley-Hamilton identity of degree
$n^{k}$ was constructed in [S1] (see also [MaMeSvW]). Since $E$ is Lie
nilpotent of index $k=2$, this identity for a matrix $A\in\mathrm{M}_{n}(E)$
is of degree $n^{2}$.

In [Do], Domokos considered a slightly modified version of the mentioned
identity, in which the left (as well as the right) coefficients are invariant
under the conjugate action of $\mathrm{GL}_{n}(K)$ on $\mathrm{M}_{n}(E)$. For
a $2\times2$\ matrix $A\in\mathrm{M}_{2}(E)$, the left scalar coefficients of
this Cayley-Hamilton identity are expressed as polynomials (over $K$) of the
traces $\mathrm{tr}(A)$, $\mathrm{tr}(A^{2})$ and $\mathrm{tr}(A^{3})$.

If $\frac{1}{2}\in R$ and $R$ satisfies the so-called weak Lie solvability
identity%
\[
\lbrack\lbrack x,y],[x,z]]=0,
\]
then for a $2\times2$\ matrix $A\in\mathrm{M}_{2}(R)$, a Cayley-Hamilton trace
identity (of degree $4$ in~$A$) with sandwich coefficients was exhibited in
[MeSvW]. If $R$ satisfies the identity%
\[
\lbrack x_{1},x_{2},...,x_{2^{s}}]_{\text{solv}}=0
\]
of general Lie solvability, then a recursive construction (also in [MeSvW])
gives a similar Cayley-Hamilton trace identity (the degree of which depends on
$s$) for a matrix $A\in\mathrm{M}_{2}(R).$

In the present paper we consider an $n\times n$ matrix $A\in\mathrm{M}_{n}(R)$
over a ring $R$ (with $1$) satisfying the identity
\[
\lbrack\lbrack x_{1},y_{1}],z_{1}][[x_{2},y_{2}],z_{2}]\cdots\lbrack\lbrack
x_{t},y_{t}],z_{t}]=0,
\]
and we prove that an invariant "power" Cayley-Hamilton identity of the form%
\[
\left(  I_{n}\lambda_{0}^{(2)}+A\lambda_{1}^{(2)}+\cdots+A^{n^{2}-1}%
\lambda_{n^{2}-1}^{(2)}+A^{n^{2}}\lambda_{n^{2}}^{(2)}\right)  ^{t}=0
\]
holds, with certain right coefficients
\[
\lambda_{i}^{(2)}\in R,\ 0\leq i\leq n^{2}-1,\quad\text{and}\quad
\lambda_{n^{2}}^{(2)}=n\bigl\{(n-1)!\bigr\}^{1+n},
\]
which are only partially determined by $A$. The cosets $\lambda_{i}^{(2)}+D,$
with $D$ the double commutator ideal $R[[R,R],R]R$ of $R$, appear in the
second right characteristic polynomial $p_{\overline{A},2}(x)$ of the natural
image $\overline{A}\in\mathrm{M}_{n}(R/D)$ of $A$ over the factor ring~$R/D$:
\[
p_{\overline{A},2}(x)=(\lambda_{0}^{(2)}+D)+(\lambda_{1}^{(2)}+D)x+\cdots
+(\lambda_{n^{2}-1}^{(2)}+D)x^{n^{2}-1}+(\lambda_{n^{2}}^{(2)}+D)x^{n^{2}}.
\]

We note that $[[x_{1},y_{1}],z_{1}][[x_{2},y_{2}],z_{2}]\cdots\lbrack\lbrack
x_{t},y_{t}],z_{t}]=0$ is a typical identity of the ring $\mathrm{U}_{t}%
(R)$\ of $t\times t$ upper triangular matrices over a ring~$R$ satisfying the
identity~$[[x,y],z]=0$ (i.e., Lie nilpotency of index 2).

Finally, using a theorem of Jennings (see [J]), we prove that if $R$\ is Lie
nilpotent of index $k$, then an identity of the form%
\[
\left(  I_{n}\lambda_{0}^{(2)}+A\lambda_{1}^{(2)}+\cdots+A^{n^{2}-1}%
\lambda_{n^{2}-1}^{(2)}+A^{n^{2}}\lambda_{n^{2}}^{(2)}\right)  ^{2^{k-2}%
}=0\text{ }(\ast)
\]
holds for $A\in\mathrm{M}_{n}(R)$. The total degree of this identity (in $A$)
is $n^{2}2^{k-2}$, a much smaller integer than the degree $n^{k}$ of $A$\ in
the right Cayley-Hamilton identity%
\[
I_{n}\lambda_{0}^{(k)}+A\lambda_{1}^{(k)}+\cdots+A^{n^{k}-1}\lambda_{n^{k}%
-1}^{(k)}+A^{n^{k}}\lambda_{n^{k}}^{(k)}=0\text{ }(\ast\ast)
\]
arising from the $k$-th right characteristic polynomial%
\[
p_{A,k}(x)=\lambda_{0}^{(k)}+\lambda_{1}^{(k)}x+\cdots+\lambda_{n^{k}-1}%
^{(k)}x^{n^{k}-1}+\lambda_{n^{k}}^{(k)}x^{n^{k}}\in R[x]
\]
of $A$ (see [S1] and [SvW1]). The advantage of $(\ast\ast)$ is that all the
coefficients are on the right side of the powers of $A$, while the expansion
of the power in $(\ast)$ yields a sum of products of the form $A^{i_{1}%
}\lambda_{i_{1}}A^{i_{2}}\lambda_{i_{2}}\cdots A^{i_{s}}\lambda_{i_{s}}$, with
$s=2^{k-2}$.

In order to provide a self-contained treatment, we present the necessary
prerequisites in sections (2) and (3).

\bigskip

\noindent2. SOME\ RESULTS\ ON\ LIE\ NILPOTENT\ RINGS

\bigskip

Let $R$ be a ring, and let $[x,y]=xy-yx$ denote the additive commutator of the
elements $x,y\in R$. It is well-known that $(R,+,[$ $,$ $])$ is a Lie ring,
$[y,x]=-[x,y]$ and $[[x,y],z]+[[y,z],x]+[[z,x],y]=0$ (the Jacobian identity).

\noindent For a sequence $x_{1},x_{2},\ldots,x_{k}$ of elements in $R$ we use
the notation $[x_{1},x_{2},\ldots,x_{k}]_{k}$\ for the left normed commutator
(Lie-)product:%
\[
\lbrack x_{1}]_{1}=x_{1}\qquad\text{ and }\qquad\lbrack x_{1},x_{2}%
,\ldots,x_{k}]_{k}=[\ldots\lbrack\lbrack x_{1},x_{2}],x_{3}],\ldots,x_{k}].
\]
Clearly, we have%
\[
\lbrack x_{1},x_{2},\ldots,x_{k},x_{k+1}]_{k+1}=[[x_{1},x_{2},\ldots
,x_{k}]_{k},x_{k+1}]=[[x_{1},x_{2}],x_{3},\ldots,x_{k},x_{k+1}]_{k}.
\]
A ring $R$ is called Lie nilpotent of index $k$ (or having property
$\mathrm{L}_{k}$) if%
\[
\lbrack x_{1},x_{2},\ldots,x_{k},x_{k+1}]_{k+1}=0
\]
is a polynomial identity on $R$. If $R$ has property $\mathrm{L}_{k}$, then
$[x_{1},x_{2},\ldots,x_{k}]_{k}$ is central for all $x_{1},x_{2},\ldots
,x_{k}\in R$.

\bigskip

A concise proof of Theorem 2.1 due to Jennings can be found in [SvW2].

\bigskip

\noindent\textbf{2.1. Theorem ([J]).} \textit{Let }$k\geq3$\textit{ be an
integer and }$R$\textit{ be a ring with }$\mathrm{L}_{k}$\textit{. Then}%
\[
\lbrack x_{1},x_{2},\ldots,x_{k}]_{k}\cdot\lbrack y_{1},y_{2},\ldots
,y_{k}]_{k}=0
\]
\textit{for all }$x_{i},y_{i}\in R$\textit{, }$1\leq i\leq k$\textit{. Thus
the two-sided ideal}%
\[
N\!=\!R\bigl\{[x_{1},x_{2},\ldots,x_{k}]_{k}\!\mid\!x_{i}\in R,1\leq i\leq
k\bigr\}\!=\!\bigl\{[x_{1},x_{2},\ldots,x_{k}]_{k}\!\mid\!x_{i}\in R,1\leq
i\leq k\bigr\}R
\]
\textit{generated by the (central) elements }$[x_{1},x_{2},\ldots,x_{k}]_{k}%
$\textit{ is nilpotent, with }\hfill\break$N^{2}=\{0\}$\textit{.}

\bigskip

\noindent\textbf{2.2. Corollary ([J]).} \textit{If }$R$\textit{ is a ring with
}$\mathrm{L}_{k}$\textit{ (}$k\geq2$\textit{), then the double commutator
ideal}%
\[
D=R[[R,R],R]R=\{%
{\textstyle\sum\nolimits_{1\leq i\leq m}}
r_{i}[[a_{i},b_{i}],c_{i}]s_{i}\mid r_{i},a_{i},b_{i},c_{i},s_{i}\in R,1\leq
i\leq m\}\vartriangleleft R
\]
\textit{is nilpotent, with }$D^{2^{k-2}}=\{0\}$\textit{.}

\bigskip

\noindent\textbf{Proof.} This follows from Theorem 2.1 by an easy induction on
$k$. $\square$

\bigskip

\noindent3. THE\ LIE NILPOTENT\ CAYLEY-HAMILTON THEOREM

\bigskip

A Lie nilpotent analogue of classical determinant theory was developed in
[S1]; further details can be found in [Do], [S2] and [SvW1]. Here we present
the basic definitions and results about the sequences of right determinants
and right characteristic polynomials, including the so-called Lie nilpotent
right Cayley-Hamilton identities.

Let $R$\ be an arbitrary (possibly non-commutative) ring or algebra with $1$,
and let
\[
\mathrm{S}_{n}=\mathrm{Sym}(\{1,\ldots,n\})
\]
denote the symmetric group of all permutations of the set $\{1,2,\ldots,n\}$.
If $A=[a_{i,j}]$ is an $n\times n$ matrix over $R$, then the element%
\[
\mathrm{sdet}(A)=\underset{\tau,\rho\in\mathrm{S}_{n}}{\sum}\mathrm{sgn}%
(\rho)a_{\tau(1),\rho(\tau(1))}\cdots a_{\tau(t),\rho(\tau(t))}\cdots
a_{\tau(n),\rho(\tau(n))}%
\]%
\[
\qquad\qquad=\underset{\alpha,\beta\in\mathrm{S}_{n}}{\sum}\mathrm{sgn}%
(\alpha)\mathrm{sgn}(\beta)a_{\alpha(1),\beta(1)}\cdots a_{\alpha(t),\beta
(t)}\cdots a_{\alpha(n),\beta(n)}%
\]
of $R$ is called the symmetric determinant of $A$.

The $(r,s)$-entry of the symmetric adjoint matrix $A^{\ast}=[a_{r,s}^{\ast}]$
of $A$ is defined as follows:%
\[
a_{r,s}^{\ast}=\underset{\tau,\rho}{\sum}\mathrm{sgn}(\rho)a_{\tau
(1),\rho(\tau(1))}\cdots a_{\tau(s-1),\rho(\tau(s-1))}a_{\tau(s+1),\rho
(\tau(s+1))}\cdots a_{\tau(n),\rho(\tau(n))}%
\]%
\[
\qquad=\underset{\alpha,\beta}{\sum}\mathrm{sgn}(\alpha)\mathrm{sgn}%
(\beta)a_{\alpha(1),\beta(1)}\cdots a_{\alpha(s-1),\beta(s-1)}a_{\alpha
(s+1),\beta(s+1)}\cdots a_{\alpha(n),\beta(n)}\text{ },
\]
where the first sum is taken over all $\tau,\rho\in\mathrm{S}_{n}$ with
$\tau(s)=s$ and $\rho(s)=r$, while the second sum is taken over all
$\alpha,\beta\in\mathrm{S}_{n}$ with $\alpha(s)=s$ and $\beta(s)=r$. We note
that the $(r,s)$ entry of $A^{\ast}$ is exactly the signed symmetric
determinant $(-1)^{r+s}\mathrm{sdet}(A_{s,r})$\ of the $(n-1)\times
(n-1)$\ minor $A_{s,r}$\ of $A$ arising from the deletion of the $s$-th row
and the $r$-th column of $A$.

The trace $\mathrm{tr}(M)$ of a matrix $M\in\mathrm{M}_{n}(R)$\ is the sum of
the diagonal entries of $M$. In spite of the fact that the well known identity
$\mathrm{tr}(AB)=\mathrm{tr}(BA)$ is no longer valid for matrices
$A,B\in\mathrm{M}_{n}(R)$ over a non-commutative $R$, we still have (see
[SvW1])%
\[
\mathrm{sdet}(A)=\mathrm{tr}(AA^{\ast})=\mathrm{tr}(A^{\ast}A).
\]
If $R$\ is commutative, then $\mathrm{sdet}(A)=n!\mathrm{\det}(A)$ and
$A^{\ast}=(n-1)!\mathrm{adj}(A)$, where $\mathrm{\det}(A)$ and $\mathrm{adj}%
(A)$ denote the ordinary determinant and adjoint, respectively, of $A$.

The right adjoint sequence $(P_{k})_{k\geq1}$ of a matrix $A\in\mathrm{M}%
_{n}(R)$ is defined by the following recursion:
\[
P_{1}=A^{\ast}\qquad\mathrm{and} \qquad P_{k+1}=(AP_{1}\cdots P_{k})^{\ast}
\]
for $k\geq1.$ The $k$-th right adjoint of~$A$ is defined as%
\[
\mathrm{radj}_{(k)}(A)=nP_{1}\cdots P_{k}.
\]
The $k$-th right determinant of $A$ is the trace of $AP_{1}\cdots P_{k}$:%
\[
\mathrm{rdet}_{(k)}(A)=\mathrm{tr}(AP_{1}\cdots P_{k}).
\]
The following theorem shows that $\mathrm{radj}_{(k)}(A)$ and $\mathrm{rdet}%
_{(k)}(A)$ can play a role similar to that played by the ordinary adjoint and
determinant, respectively, in the commutative case.

\bigskip

\noindent\textbf{3.1. Theorem ([S1], [SvW1]).}\textit{ If }$\frac{1}{n}\in
R$\textit{ and the ring }$R$\textit{ is Lie nilpotent of index }~$k$\textit{,
then for a matrix }$A\in\mathrm{M}_{n}(R)$\textit{\ we have}%
\[
A\mathrm{radj}_{(k)}(A)=nAP_{1}\cdots P_{k}=\mathrm{rdet}_{(k)}(A)I_{n}.
\]

\bigskip

The above Theorem 3.1 is not used explicitly in the sequel, however it helps
our understanding and serves as a starting point in the proof of Theorem 3.3.

Let $R[x]$ denote the ring of polynomials in the single commuting
indeterminate~$x$, with coefficients in $R$. The $k$-th right characteristic
polynomial of $A$ is the $k$-th right determinant of the $n\times n$ matrix
$xI_{n}-A$ in $\mathrm{M}_{n}(R[x])$:%
\[
p_{A,k}(x)=\mathrm{rdet}_{(k)}(xI_{n}-A).
\]

\bigskip

\noindent\textbf{3.2. Proposition ([S1], [SvW1]).}\textit{ The }$k$\textit{-th
right characteristic polynomial}

\noindent$p_{A,k}(x)\in R[x]$\textit{ of }$A\in\mathrm{M}_{n}(R)$\textit{\ is
of the form}%
\[
p_{A,k}(x)=\lambda_{0}^{(k)}+\lambda_{1}^{(k)}x+\cdots+\lambda_{n^{k}-1}%
^{(k)}x^{n^{k}-1}+\lambda_{n^{k}}^{(k)}x^{n^{k}},
\]
\textit{where }$\lambda_{0}^{(k)},\lambda_{1}^{(k)},\ldots,\lambda_{n^{k}%
-1}^{(k)},\lambda_{n^{k}}^{(k)}\in R$\textit{ and }$\lambda_{n^{k}}%
^{(k)}=n\bigl\{  (n-1)!\bigr\}  ^{1+n+n^{2}+\cdots+n^{k-1}}$\textit{.}

\bigskip

\noindent\textbf{3.3. Theorem ([S1], [SvW1]).}\textit{ If }$\frac{1}{n}\in
R$\textit{ and the ring }$R$\textit{ is Lie nilpotent of index }~$k$\textit{,
then a right Cayley-Hamilton identity}%
\[
(A)p_{A,k}=I_{n}\lambda_{0}^{(k)}+A\lambda_{1}^{(k)}+\cdots+A^{n^{k}-1}%
\lambda_{n^{k}-1}^{(k)}+A^{n^{k}}\lambda_{n^{k}}^{(k)}=0
\]
\textit{with right scalar coefficients holds for }$A\in\mathrm{M}_{n}%
(R)$\textit{. We also have }$(A)u=0$\textit{, where }$u(x)=p_{A,k}%
(x)h(x)$\textit{ and }$h(x)\in R[x]$\textit{ is arbitrary.}

\bigskip

\noindent\textbf{3.4. Theorem ([Do]). }\textit{If }$\frac{1}{2}\in R$\textit{
and the ring }$R$\textit{ is Lie nilpotent of index }$2$\textit{, then for a
}$2\times2$\textit{ matrix }$A\in\mathrm{M}_{2}(R)$\textit{\ the right
Cayley-Hamilton identity in the above 3.3 can be written in the following
trace form:}%
\[
(A)p_{A,2}\!=\!\!I_{2}\!\!\left(  \!\frac{1}{2}\mathrm{tr}^{4}(A)\!+\!\frac
{1}{2}\mathrm{tr}^{2}(A^{2})\!+\!\frac{1}{4}\mathrm{tr}^{2}(A)\mathrm{tr}%
(A^{2})\!-\!\frac{5}{4}\mathrm{tr}(A^{2})\mathrm{tr}^{2}%
(A)\!+\!\bigl[\mathrm{tr}(A^{3}),\mathrm{tr}(A)\bigr]\!\right)  \!+
\]%
\[
A\Bigl(\!\mathrm{tr}(A)\mathrm{tr}(A^{2})+\mathrm{tr}(A^{2})\mathrm{tr}%
(A)-2\mathrm{tr}^{3}(A)\!\Bigr)+A^{2}\Bigl(\!4\mathrm{tr}^{2}(A)-2\mathrm{tr}%
(A^{2})\!\Bigr)-A^{3}\Bigl(\!4\mathrm{tr}(A)\!\Bigr)+2A^{4}=0.
\]

\bigskip

\noindent\textbf{3.5. Corollary ([Do]).}\textit{ If }$\frac{1}{2}\in
R$\textit{ and the ring }$R$\textit{ is Lie nilpotent of index }~$2$\textit{,
then, for every }$A\in\mathrm{M}_{2}(R)$,
\[
\mathrm{tr}(A)=\mathrm{tr}(A^{2})=0\quad\text{\textit{implies that}}\quad
A^{4}=0\mathit{.}%
\]

\newpage

\noindent4. MATRICES OVER A\ RING\ WITH $[[x_{1},y_{1}],z_{1}][[x_{2}%
,y_{2}],z_{2}]\cdots\lbrack\lbrack x_{t},y_{t}],z_{t}]=0$

\bigskip

We shall make use of the following well-known fact.

\bigskip

\noindent\textbf{4.1. Proposition.}\textit{ If }$[[x_{1},y_{1}],z_{1}%
][[x_{2},y_{2}],z_{2}]\cdots\lbrack\lbrack x_{t},y_{t}],z_{t}]=0$\textit{ is a
polynomial identity on a ring }$R$\textit{, then }$D^{t}=\{0\}$\textit{, with
}$D$\textit{ the ideal }$R[[R,R],R]R$\textit{ of }$R$\textit{.}

\bigskip

\noindent\textbf{4.2. Theorem.}\textit{ If }$\frac{1}{2}\in R$ \textit{and
}$A\in\mathrm{M}_{n}(R)$\textit{ is a matrix over a ring }$R$\textit{
satisfying the polynomial identity }$[[x_{1},y_{1}],z_{1}][[x_{2},y_{2}%
],z_{2}]\cdots\lbrack\lbrack x_{t},y_{t}],z_{t}]=0$\textit{, then an invariant
"power" Cayley-Hamilton identity of the form}%
\[
\left(  I_{n}\lambda_{0}^{(2)}+A\lambda_{1}^{(2)}+\cdots+A^{n^{2}-1}%
\lambda_{n^{2}-1}^{(2)}+A^{n^{2}}\lambda_{n^{2}}^{(2)}\right)  ^{t}=0
\]
\textit{holds, with certain right coefficients }
\[
\lambda_{i}^{(2)}\in R,\ 0\leq i\leq n^{2}-1,\quad and\quad\lambda_{n^{2}%
}^{(2)}=n\bigl\{(n-1)!\bigr\}^{1+n}%
\]
\textit{(only partially determined by }$A$\textit{). The cosets }$\lambda
_{i}^{(2)}+D$\textit{ with }$D=R[[R,R],R]R\vartriangleleft R$\textit{ appear
in the second right characteristic polynomial }$p_{\overline{A},2}(x)$\textit{
of the natural image}

\noindent$\overline{A}\in M_{n}(R/D)$\textit{ of }$A$\textit{ over the factor
ring }$R/D$\textit{:}%
\[
p_{\overline{A},2}(x)=(\lambda_{0}^{(2)}+D)+(\lambda_{1}^{(2)}+D)x+\cdots
+(\lambda_{n^{2}-1}^{(2)}+D)x^{n^{2}-1}+(\lambda_{n^{2}}^{(2)}+D)x^{n^{2}}%
\in(R/D)[x].
\]

\bigskip

\noindent\textbf{Proof.} Consider the factor ring $R/D$, where
$D=R[[R,R],R]R\vartriangleleft R$ is the double commutator ideal. If
$A=[a_{i,j}]\in\mathrm{M}_{n}(R)$, then we use the notation $\overline
{A}=[a_{i,j}+D]$ for the image of $A$ in $\mathrm{M}_{n}(R/D)$. Since $R/D$ is
Lie nilpotent of index $2$, Theorem~3.3 implies that, in $\mathrm{M}_{n}%
(R/D)$,%
\[
(\overline{A})p_{\overline{A},2}=\overline{I_{n}}(\lambda_{0}^{(2)}%
+D)+\overline{A}(\lambda_{1}^{(2)}+D)+\cdots+(\overline{A})^{n^{2}-1}%
(\lambda_{n^{2}-1}^{(2)}+D)+(\overline{A})^{n^{2}}(\lambda_{n^{2}}%
^{(2)}+D)=\overline{0},
\]
where%
\[
p_{\overline{A},2}(x)=\mathrm{rdet}_{(k)}(x\overline{I_{n}}-\overline{A})
\]%
\[
=(\lambda_{0}^{(2)}+D)+(\lambda_{1}^{(2)}+D)x+\cdots+(\lambda_{n^{2}-1}%
^{(2)}+D)x^{n^{2}-1}+(\lambda_{n^{2}}^{(2)}+D)x^{n^{2}}%
\]
is the second right characteristic polynomial of $\overline{A}$ in $(R/D)[x]$.
Clearly,%
\[
\overline{I_{n}\lambda_{0}^{(2)}+A\lambda_{1}^{(2)}+\cdots+A^{n^{2}-1}%
\lambda_{n^{2}-1}^{(2)}+A^{n^{2}}\lambda_{n^{2}}^{(2)}}%
\]%
\[
=\overline{I_{n}}(\lambda_{0}^{(2)}+D)+\overline{A}(\lambda_{1}^{(2)}%
+D)+\cdots+(\overline{A})^{n^{2}-1}(\lambda_{n^{2}-1}^{(2)}+D)+(\overline
{A})^{n^{2}}(\lambda_{n^{2}}^{(2)}+D)=\overline{0}%
\]
implies that%
\[
I_{n}\lambda_{0}^{(2)}+A\lambda_{1}^{(2)}+\cdots+A^{n^{2}-1}\lambda_{n^{2}%
-1}^{(2)}+A^{n^{2}}\lambda_{n^{2}}^{(2)}\in\mathrm{M}_{n}(D).
\]
Now $D^{t}=\{0\}$ is a consequence of Proposition 4.1, whence $\left(
\mathrm{M}_{n}(D)\right)  ^{t}=\{0\}$ and%
\[
\left(  I_{n}\lambda_{0}^{(2)}+A\lambda_{1}^{(2)}+\cdots+A^{n^{2}-1}%
\lambda_{n^{2}-1}^{(2)}+A^{n^{2}}\lambda_{n^{2}}^{(2)}\right)  ^{t}=0
\]
follows. $\square$

\bigskip

\noindent\textbf{4.3. Remark.} If $[x_{1},y_{1}][x_{2},y_{2}]\cdots\lbrack
x_{t},y_{t}]=0$ is a polynomial identity on $R$ and $A\in\mathrm{M}_{n}(R)$,
then using the commutator ideal $T=R[R,R]R$ and the natural image
$\widetilde{A}\in\mathrm{M}_{n}(R/T)$\ of $A$ over the commutative ring $R/T$,
a similar argument as in the proof of Theorem 4.2 gives that%
\[
\left(  I_{n}\lambda_{0}^{(1)}+A\lambda_{1}^{(1)}+\cdots+A^{n-1}\lambda
_{n-1}^{(1)}+A^{n}\lambda_{n}^{(1)}\right)  ^{t}=0
\]
holds, where $p_{\widetilde{A},1}(x)=(\lambda_{0}^{(1)}+T)+(\lambda_{1}%
^{(1)}+T)x+\cdots+(\lambda_{n-1}^{(1)}+T)x^{n-1}+(\lambda_{n}^{(1)}+T)x^{n}$
is the $n!$ times scalar multiple of the classical characteristic polynomial
of $\widetilde{A}$ in $(R/T)[x]$ with $\lambda_{n}^{(1)}=n!$.

\bigskip

\noindent\textbf{4.4. Theorem.}\textit{ If }$\frac{1}{2}\in R$ \textit{and
}$A\in\mathrm{M}_{n}(R)$\textit{ is a matrix over a Lie nilpotent ring }%
$R$\textit{ of index }$k$\textit{, then an invariant "power" Cayley-Hamilton
identity of the form}%
\[
\left(  I_{n}\lambda_{0}^{(2)}+A\lambda_{1}^{(2)}+\cdots+A^{n^{2}-1}%
\lambda_{n^{2}-1}^{(2)}+A^{n^{2}}\lambda_{n^{2}}^{(2)}\right)  ^{2^{k-2}}=0
\]
\textit{holds, with certain right coefficients}%
\[
\lambda_{i}^{(2)}\in R,\ 0\leq i\leq n^{2}-1,\quad\text{\textit{and}}%
\quad\lambda_{n^{2}}^{(2)}=n\bigl\{(n-1)!\bigr\}^{1+n}%
\]
\textit{(only partially determined by }$A$\textit{). The cosets }$\lambda
_{i}^{(2)}+D$\textit{ with }$D=R[[R,R],R]R\vartriangleleft R$\textit{ appear
in the second right characteristic polynomial }$p_{\overline{A},2}(x)$\textit{
of the natural image}

\noindent$\overline{A}\in M_{n}(R/D)$\textit{ of }$A$\textit{ over the factor
ring }$R/D$\textit{:}%
\[
p_{\overline{A},2}(x)=(\lambda_{0}^{(2)}+D)+(\lambda_{1}^{(2)}+D)x+\cdots
+(\lambda_{n^{2}-1}^{(2)}+D)x^{n^{2}-1}+(\lambda_{n^{2}}^{(2)}+D)x^{n^{2}}%
\in(R/D)[x].
\]

\bigskip

\noindent\textbf{Proof.} According to Jennings's result (Corollary 2.2), the
double commutator ideal%
\[
D=R[[R,R],R]R=\bigl\{%
{\textstyle\sum\nolimits_{1\leq i\leq m}}
r_{i}[[a_{i},b_{i}],c_{i}]s_{i}\mid r_{i},a_{i},b_{i},c_{i},s_{i}\in R,1\leq
i\leq m\bigr\}\vartriangleleft R
\]
is nilpotent, with $D^{2^{k-2}}=\{0\}$. Thus the application of Theorem 4.2
gives our identity. $\square$

\bigskip

\noindent\textbf{4.5. Remark.} If $k=2$, then $R[[R,R],R]R=\{0\}$, and the
identity in Theorem 4.4 remains the same as the Lie nilpotent right
Cayley-Hamilton identity in Theorem~3.3.

\bigskip

\noindent\textbf{4.6. Remark.} The Grassmann algebra
\[
E=K\left\langle v_{1},v_{2},...,v_{i},...\mid v_{i}v_{j}+v_{j}v_{i}=0\text{
for all }1\leq i\leq j\right\rangle
\]
over a field $K$ (with $2\neq0$) has property $\mathrm{L}_{2}$, and%
\[
\lbrack v_{1},v_{2}]\cdot\lbrack v_{3},v_{4}]\cdot\ \cdots\ \cdot\lbrack
v_{2t-1},v_{2t}]=2^{t}v_{1}v_{2}\cdots v_{2t}\neq0
\]
shows that $\mathrm{L}_{2}$ does not imply the identity $[x_{1},y_{1}%
][x_{2},y_{2}]\cdots\lbrack x_{t},y_{t}]=0$ for any $t$. Thus the identity
mentioned in Remark 4.3 cannot be used directly to derive new identities for
matrices over a Lie nilpotent ring of index $k\geq2$. However, as the referee
pointed out, the following (weak) version of Latyshev's theorem provides a
possibility to use Remark 4.3 in order to obtain another remarkable "power"
Cayley-Hamilton identity.

\bigskip

\textit{Theorem ([L]). If }$S$\textit{ is a Lie nilpotent algebra (over an
infinite field) of index }$k$\textit{, generated by }$m$\textit{ elements,
then there exists an integer }$d=d(k,m)$\textit{ such that }$S$%
\textit{\ satisfies the polynomial identity }$[x_{1},y_{1}][x_{2},y_{2}%
]\cdots\lbrack x_{d},y_{d}]=0$\textit{. (In the original version }%
$S$\textit{\ satisfies a so-called nonmatrix polynomial identity.)}

\bigskip

\noindent If $A\in\mathrm{M}_{n}(R)$ is a matrix over a Lie nilpotent algebra
(over an infinite field) $R$ of index $k$, then $A\in\mathrm{M}_{n}(S)$, where
$S$ is the (unitary) subalgebra generated by the $n^{2}$ entries of $A$. Thus
$[x_{1},y_{1}][x_{2},y_{2}]\cdots\lbrack x_{d},y_{d}]=0$ is a polynomial
identity on $S$ with $d=d(k,n^{2})$ and Remark 4.3 gives that an identity
\[
\left(  I_{n}\lambda_{0}^{(1)}+A\lambda_{1}^{(1)}+\cdots+A^{n-1}\lambda
_{n-1}^{(1)}+A^{n}\lambda_{n}^{(1)}\right)  ^{d(k,n^{2})}=0
\]
of degree $nd(k,n^{2})$\ holds. Unfortunately, our knowledge about
$d(k,n^{2})$\ is very limited, the fact that $d(2,4)=3$\ was mentioned by the referee.

\bigskip

\noindent\textbf{4.7. Remark.} If $R$ is an algebra over a field $K$\ of
characteristic zero, then the invariance of the identities in 4.2 and 4.4
means that $p_{\overline{T^{-1}AT},2}(x)=p_{\overline{A},2}(x)$ holds for any
$T\in\mathrm{GL}_{n}(K)$ (see [Do]).

\bigskip

\noindent\textbf{4.8. Corollary.}\textit{ If }$\frac{1}{2}\in R$\textit{ and
the ring }$R$\textit{ is Lie nilpotent of index }~$k$\textit{, then, for every
}$A\in\mathrm{M}_{2}(R)$,
\[
\mathrm{tr}(A)=\mathrm{tr}(A^{2})=0\quad\text{\textit{implies that}}\quad
A^{2^{k}}=0\mathit{.}%
\]

\bigskip

\noindent\textbf{Proof.} Using $D=R[[R,R],R]R\vartriangleleft R$,
$\overline{A}\in\mathrm{M}_{2}(R/D)$ and%
\[
\mathrm{tr}(\overline{A})=\mathrm{tr}(A)+D=0,\text{ }\mathrm{tr}((\overline
{A})^{2})=\mathrm{tr}(\overline{A^{2}})=\mathrm{tr}(A^{2})+D=0,
\]
the application of Corollary 3.5\ ensures that $\overline{A^{4}}=\left(
\overline{A}\right)  ^{4}=\overline{0}$. Thus the nilpotency of $D$
($D^{2^{k-2}}=\{0\}$) gives that $A^{2^{k}}=\left(  A^{4}\right)  ^{2^{k-2}%
}=0$. $\square$

\bigskip

\noindent\textbf{4.9. Remark.} According to the following important
observation of the referee, the use of Latyshev's theorem gives an $n\times n$
variant of Corollary 4.8. If $A\in\mathrm{M}_{n}(R)$ is a matrix over a Lie
nilpotent algebra (over a field $K$\ of characteristic zero) $R$ of index $k$,
then we prove that%
\[
\mathrm{tr}(A)=\mathrm{tr}(A^{2})=\cdots=\mathrm{tr}(A^{n})=0
\]
implies that $A^{nd(k,n^{2})}=0$. Indeed, $A\in\mathrm{M}_{n}(S)$, where
$S\subseteq R$ is the (unitary) subalgebra of $R$ generated by the entries of
$A$. Now consider the natural image $\widetilde{A}\in\mathrm{M}_{n}%
(S/S[S,S]S)$ of $A$. The application of the well-known fact that%
\[
\mathrm{tr}(\widetilde{A})=\mathrm{tr}((\widetilde{A})^{2})=\cdots
=\mathrm{tr}((\widetilde{A})^{n})=\widetilde{0}%
\]
implies that $\widetilde{A^{n}}=(\widetilde{A})^{n}=\widetilde{0}$ (it is a
consequence of the Newton trace formulae for the coefficients of the
characteristic polynomial $p_{\widetilde{A},1}(x)\in(S/S[S,S]S)[x]$, where the
factor $S/S[S,S]S$ is a commutative algebra over $K$). Since
$(S[S,S]S)^{d(k,n^{2})}=\{0\}$ by Latyshev's theorem and $A^{n}\in
\mathrm{M}_{n}(S[S,S]S)$, we obtain the desired equality.

\bigskip

\noindent ACKNOWLEDGMENT

\noindent The authors thank the referee for the valuable report and the
important contributions mentioned in Remarks 4.6 and 4.9.

\bigskip

\noindent REFERENCES

\bigskip

\noindent\lbrack Do] M. Domokos, \textit{Cayley-Hamilton theorem for }%
$2\times2$\textit{\ matrices over the Grassmann algebra}, J. Pure
Appl.~Algebra 133 (1998), 69-81.

\noindent\lbrack Dr] V. Drensky, \textit{Free algebras and PI-algebras}.
Graduate course in algebra. Springer-Verlag Singapore, Singapore, 2000.

\noindent\lbrack DrF] V. Drensky and E. Formanek, \textit{Polynomial identity
rings}. Advanced Courses in Mathematics. CRM Barcelona. Birkh\"{a}user-Verlag,
Basel, 2004.

\noindent\lbrack J] S. A. Jennings, \textit{On rings whose associated Lie
rings are nilpotent}, Bull. Amer. Math. Soc. 53 (1947), 593-597.

\noindent\lbrack K] A. R. Kemer, \textit{Ideals of identities of associative
algebras}. Translated from the Russian by C. W. Kohls. Translations of
Math.~Monographs, 87. American Mathematical Society, Providence, RI, 1991.

\noindent\lbrack L] V. N. Latyshev,\textit{ Generalization of the Hilbert
theorem on the finiteness of bases }(Russian), Sib. Mat. Zhurn. 7 (1966),
1422-1424. Translation: Sib. Math. J. 7 (1966), 1112-1113.

\noindent\lbrack MaMeSvW] L. M\'arki, J. Meyer, J. Szigeti and L. van Wyk,
\textit{Matrix representations of finitely generated Grassmann algebras and
some consequences}, Israel J.~Math.~208 (2015), 373-384.

\noindent\lbrack MeSvW] J. Meyer, J. Szigeti and L. van Wyk, \textit{A
Cayley-Hamilton trace identity for }$2\times2$\textit{ matrices over
Lie-solvable rings}, Linear Algebra Appl.~436 (2012), 2578-2582.

\noindent\lbrack PaSch] R. Par\'{e} and W. Schelter, \textit{Finite extensions
are integral}, J.~Algebra 53 (1978), 477-479.

\noindent\lbrack Pe1] K. R. Pearson, \textit{A lower bound for the degree of
polynomials satisfied by matrices}, J.~Aust.~Math.~Soc.~Ser.~A 27 (1979), 430-436.

\noindent\lbrack Pe2] K. R. Pearson, \textit{Degree 7 monic polynomials
satisfied by a }$3\times3$\textit{ matrix over a noncommutative ring},
Comm.~Algebra 10 (1982), 2043-2073.

\noindent\lbrack Rob] J. C. Robson, \textit{Polynomials satisfied by
matrices}, J.~Algebra 55 (1978), 509-520.

\noindent\lbrack Row] L. H. Rowen, \textit{Polynomial identities in ring
theory}. Pure and Applied Mathematics, 84. Academic Press, New York-London, 1980.

\noindent\lbrack S1] J. Szigeti, \textit{New determinants and the
Cayley-Hamilton theorem for matrices over Lie nilpotent rings},
Proc.~Amer.~Math.~Soc.~125 (1997), 2245-2254.

\noindent\lbrack S2] J. Szigeti, \textit{On the characteristic polynomial of
supermatrices}, Israel J.~Math.~107 (1998), 229-235.

\noindent\lbrack SvW1] J. Szigeti and L. van Wyk,\textit{ Determinants for
}$n\times n$\textit{ matrices and the symmetric Newton formula in the
}$3\times3$\textit{ case,} Linear Multilinear Algebra~ 62 (2014), 1076-1090.

\noindent\lbrack SvW2] J. Szigeti and L. van Wyk,\textit{ On Lie nilpotent
rings and Cohen's theorem,} Comm.~Algebra 43 (2015), 4783--4796.

\end{document}